\documentclass[a4paper,12pt]{article}
\usepackage{graphicx,amssymb}
\usepackage{amsmath}
\usepackage{amssymb}
\usepackage{amsthm}
\usepackage{amsfonts}
\usepackage{graphicx}
\newtheorem{corollary}{Corollary}
\newtheorem{theorem}[corollary]{Theorem}
\newtheorem{lemma}[corollary]{Lemma}

\newtheorem{proposition}[corollary]{Proposition}

\begin{document}

\font \m=msbm10
\newcommand{\R}{{\hbox {\m R}}}
\newcommand{\C}{{\hbox {\m C}}}
\newcommand{\Z}{{\hbox {\m Z}}}
\newcommand{\N}{{\hbox {\m N}}}
\newcommand{\U}{{\hbox {\m U}}}
\newcommand{\HH}{{\hbox {\m H}}}
\newcommand{\aba}{{\alpha}}

\title{Some locally self-interacting walks on the integers}
\author{Anna Erschler\footnote {CNRS, Universit\'e Paris-Sud 11}
 \and B\'alint T\'oth\footnote {Technical University Budapest}
 \and Wendelin Werner\footnote {Universit\'e Paris-Sud 11 and Ecole Normale Sup\'erieure}}
\date {}
\maketitle
\def \interior {{\hbox {int}}}
\def \eps {\varepsilon}

\begin {abstract}
We study certain self-interacting walks on the set of integers, that choose to jump to the right or to the left randomly but influenced by the number of times they have
previously jumped along the edges in the finite neighbourhood of their current position (in the present paper, typically, we will discuss the case where one considers the neighbouring edges and the next-to-neighbouring edges). We survey a variety of possible behaviours, including some where the walk is 
eventually confined to an interval of large length.
We also focus on certain ``asymmetric''  drifts, where we prove
that with positive probability, the walks behave deterministically on large scale and move like  $n \mapsto c \sqrt {n}$ or like $n \mapsto c \log n$.
\end {abstract}

\bigbreak
\begin {center}
{\sl Dedicated to Erwin Bolthausen on the occasion of his 65th birthday}
\end {center}

\bigbreak

\section {Introduction}

A locally self-interacting random walk on the set of integers is a sequence of integer-valued random variables $(X_n, n \ge 0)$ with $|X_{n+1} - X_n |= 1$ for all $n \ge 0$, that can be defined inductively as follows:
One initializes a ``local time profile'' by first choosing some function $L_0$ that associates to each edge $e$ of the lattice $\Z$ a real number $L_0 (e)$
(the set $E$ of unoriented edges can be identified to the set $\Z + 1/2$ or to the set of couples $\{x,x+1\}$ where $x \in \Z$) and one chooses the starting point
$X_0$. The simplest choice is of course to set this initial function $L_0$ to be equal to zero on all edges, and to start at the origin (i.e., $X_0=0$)
with probability one.

The law of the walker is then described via some function $R$ from $\R^E$ to $[0,1]$. This function describes the probability of jumping to the right in terms of  the local time profile ``seen from the walker at its current position''. More precisely, one defines inductively, for each $n \ge 0$, $X_{n+1}$ to be equal to either $X_n +1$ or $X_n -1 $ in such a way that
$$
P ( X_{n+1} =  X_n + 1 \ | \ X_0, \ldots, X_n ) = R ( (\ell_n (e), e \in E)) ,$$
 where
$$\ell_n (\cdot) = L_n ( X_n + \cdot)$$ is the local time profile seen from the particle at time $n$.
At each step, one updates this local time profile as follows:
$e_{n+1}= \{ X_n , X_{n+1} \} $ denotes the edge of the $(n+1)$-th jump and
$$ L_{n+1} ( e) = L_{n} (e) + 1_{\{e=e_{n+1}\}} $$
is the new local time profile. When $L_0$ is identically zero, then $L_n (e)$ denotes the number of times the walk has jumped on the edge $e$ before time $n$.
This procedure clearly describes completely the law of the family $(X_n , n \ge 0)$.

Similar definitions work of course on other graphs than $\Z$ (for instance on $\Z^d$).
In the present paper, we will however restrict our discussion to the one-dimensional case. One can also use local times on sites instead of edges
but the latter case turns out to be often easier to handle.

Depending on the function $R$, the walker can be of self-repelling type (it prefers to go into the direction it has less visited in the past), self-attracting, or more complicated (it can be for instance repelled by the visits to the neighbouring edge, but attracted by its past visits to further edges).

When the function $R$ is ``cylindric'' i.e. when $R( (\ell (e), e \in E))$ is a function on the values $\ell (e_1), \ldots , \ell (e_K)$ for a finite set of edges $e_1, \ldots, e_K$, the interaction can be described as  ``local''. In other words, for some constant $K'$, the past trajectory $(X_0, X_1, \ldots, X_n)$ influences the outcome of
$X_{n+1}-X_n$ only via the values of $L_n (e)$ for the edges in the $K'$-neighbourhood of $X_n$. Note that in this cylindric case, if there is a large scale limit to the process $(X_n, n \ge 0)$ its evolution will be  ``purely'' local as the finite-range becomes scaled down.

Such locally self-interacting walks provide a particular class of random evolutions of functions as the sequence of functions
 $(\ell_n(\cdot))_{n \ge 0}$ form a Markov process. It is in fact very natural to study the evolutions where the function $R$ is also height-invariant i.e. the particular case where
$$ R ((\ell (e), e \in E)) = R ((l+ \ell (e), e \in E))$$
for all $(\ell(\cdot))$ and all real $l$. In other words, the walk is sensitive to the local discrete gradient of $\ell$ but not to its actual height. This excludes for instance the family of ``linear'' reinforcements (see for instance \cite {Pe} for a survey of these other types of reinforcements). Another way to view these evolutions is to define, for all $x \in \Z$ and $n \ge 0$, the local gradient
$$ \eta_n (x) = \nabla \ell_n (x) := \ell_n (x+1/2) - \ell_n (x-1/2)$$
and to say that $R((\ell))$ is in fact a function $\tilde R$ of $(\eta)$.

In \cite {T1}, a case where $R (( \ell))$ was just a function of
$\eta (0) := \ell (1/2) - \ell (-1/2)$ was studied. Loosely speaking, this walk (sometimes referred to as true self-avoiding walk or true self-repelling walk -- TSRW) is driven by the negative gradient of the local time at its actual position i.e., the walk prefers to go in the direction it has visited less often in the past.
It turns out that a non-trivial scaling shows up, as $X_n / n^{2/3}$ converges in law to some non-trivial random variable. In fact, the entire limiting process can be described \cite {TW} (it is called the true self-repelling motion -- TSRM). It  corresponds to a random and local evolution of functions, that is neither a deterministic partial differential equation nor constructed via a white-noise driven stochastic partial differential equation.

This raises the questions whether other such processes exist, whether the previous scaling limit is ``stable'' (do all similar discrete models have the same limiting behaviour? What properties of the discrete model is actually captured by this scaling limit?), and
 provides a motivation to study self-interacting random walks associated to more general functions $R$.

The family of models that we will focus on in the present paper
are those where the self-interaction is governed by a linear combination of the values of the local time. In other words, we let
$$ D ((\ell))= \sum_{e \in E } a_e \ell (e)   $$
where all but finitely many $a_e$'s are zero,
and we choose the local ``drift'' at step $n$ according to the value of $D ((\ell_n))$. The fact that this evolution in height-invariant corresponds to the condition
$$\sum_{e \in E} a_e = 0.$$
Then, we take
$$
R (( \ell)) = F ( D ((\ell))) $$
for some increasing function $F$ such that $F ( -y) = 1 - F(y)$. In other words, if $D$ is (very) positive (respectively,  negative), then the walk will tend to jump to the right (respectively, to the left), so that $D$ indeed plays the role of a local drift.
A particularly natural choice for $F$ turns out to be
$$F(y) = \frac {e^y}{ e^y + e^{-y}}.$$

In order to start investigating the variety of possible behaviours, one can for instance start to look at the case where $R(\ell)$ is a function of
$(\ell (-3/2)$, $\ell (-1/2)$, $\ell(1/2)$, $\ell(3/2))$.
A first natural choice is to pick two real numbers $a$ and $b$, and to drive the walker by the linear combination
$$ D^{a,b} ((\ell)) = D^{a,b} ((\ell)) := a \ell (-3/2) + b \ell (-1/2) - b \ell (1/2) - a \ell (3/2).$$
 Note that these are the choices such that both height-invariance (because the sum of the terms $a+b-b-a$ is zero) and {\em left-right symmetry}
(the walks $(X_n)$ and $(-X_n)$ have the same law because the coefficients in front of $\ell(e)$ and $\ell (-e)$ are opposite) hold.

Note that when $b$ is positive and $a=0$, this is the previously mentioned TSRW  with TSRM as scaling limit.
 When $b >0$ and $|a|$ is very small, one can therefore view this process as a perturbation of the TSRW.
At
the other end of the picture, when $b$ is negative, it is easy to check that the walk will (with positive probability) eventually jump back and forth on one single edge (the more it jumps on it, the more it is attracted by it) which  is not a particularly interesting case.
In Section \ref {S2}, we will describe various possible asymptotic behaviours, depending on the values of $a$ and $b$.
We shall see that intuition can sometimes be misleading; repellance by second neighbours does not always have the same effect as repellance by immediate neighbours as it can create traps. One interesting feature that we will mention is the existence (when $b>0> a$) of a phase transition between stuck walks (that can eventually remain in a finite interval) and the non-stuck ones (with almost surely an infinite range), depending on whether $-a > b/3$ or not.

A second natural possibility is to break the left-right symmetry, and to drive the walker by a linear combination such as
$$ - \ell (-3/2) + \ell (-1/2) + \ell (1/2) -  \ell (3/2)$$
which will be among those studied in Section \ref {S3} (with deterministic transient behaviour of the walker).

Section \ref {S2} is more of survey-type presenting results proved in \cite {ETW}, as well as heuristics and conjectures, while Section \ref {S3},
will contain original results with actual proofs.

Throughout the paper, we will mention results that hold with positive probability, like ``with positive probability, the walk will eventually remain stuck on a single edge''. Clearly, one then implicitly conjectures that this is the generic behaviour of the walk (i.e. that the corresponding asymptotic event will hold almost surely), but in general, proving such almost sure results for self-interacting walks is not a straightforward task (and it requires different arguments), see for instance \cite {Ta} for the linearly reinforced walk. We will not focus on these questions in the present paper.

We would like to emphasize here that (despite the ``phase diagram'' that we will draw at the end of the next section) the aim of the present paper is not to give a full classification of the possible asymptotic behaviour of all these cylindric self-interacting walks. Our goal is rather to point out certain phenomenological features, in order to shed some light on how sensitive (or not) these walks are with respect to their microscopic definition, and what feature is kept in the scaling limit.

\medbreak
\noindent
{\bf A warm-up.}
As a warm-up to maybe help (if it does not help the reader, then she/he should not hesitate to skip this paragraph as it will not be essential later on) the
following analysis, let us make some non-rigorous analogy with Fibonacci-type sequences and their related polynomials.
Recall that our models are ``driven'' by a linear combination of the type
$$ D ((\ell)) = \sum_{k=k_0}^{k_1} a_{k + 1/2} \ell (k+1/2),$$
where $a_{k_0} \not= 0 \not= a_{k_1}$.
If a walk is in some local equilibrium state (it does not ``shoot in one direction'', it is somewhere in the middle of its past range, and the values of $\ell$ are all very large), then one could expect that the value of $D$ is rather small at most of the sites (at least compared to the total number of steps). After renormalization by the number of steps, this would indicate that
$$\sum_{k=k_0}^{k_1} a_{k+1/2} \ell_n (j+k+1/2) << \ell_n (1/2)$$ i.e. that $j \mapsto F_j = \ell_n (j+1/2) / \ell_n (1/2)$  follows some Fibonacci-type progression. As we know, the behaviour of such
progressions is controlled by the values (and position with respect to the unit circle) of the roots of the corresponding polynomial $P(x)$ where 
\begin{equation}
\label{poly}
P(x) x^{k_0} = \sum_{k \ge k_0}  a_{k+1/2} x^k.
\end{equation}
The fact that the dynamics is height-invariant corresponds to the case where $1$ is a root of $P$.
Hence, it is in a way more appropriate to describe the self-interaction by its polynomial (in the case where the self-interaction is left-right symmetric, it does indeed characterize the self-interaction), and we will use this when we will illustrate some of our simulations. Note that  the function $\theta\mapsto P(e^{i\theta})$, defined for $\theta\in[-\pi,\pi)$ is actually the Fourier transform of the sequence $(a_k)$. Thus, the qualitative behaviour of the walk will depend on the type of the Fibonacci polynomial $P(x)$, or, what is equivalent of the Fourier transform of the sequence of coefficients  $(a_k)$. E.g., positive definiteness of this sequence plays an essential role.

\section {Survey of left-right symmetric cases}
\label {S2}

During most of this section, the setup is the one described in the introduction, where the walk is ``driven'' by
$D^{a,b} ((\ell))$, i.e. when we are considering the next-to-neighbouring interactions that are both height-invariant and left-right symmetric
(we will comment on interactions with larger range at the end of the section).
The meaning of $a$ and $b$ is as follows:
\begin {itemize}
 \item If $b$ is positive (respectively,  negative), then the walk is repelled (respectively,  attracted) by its previous visits to the neighbouring edge.
\item $a$ plays the same role as $b$ but for the next-to-neigboring edges.
\item $|a|$ and $|b|$ describe the intensity of these self-interactions.
\end {itemize}
As we shall see, in most of these cases, the qualitative asymptotic behaviour of the walk will depend on the signs of $a$ and $b$, and on the ratio between $|a|$ and $|b|$.

We can divide the parameter space into several classes. Let us first notice that in the case where $b<0$, then with positive probability, the walk $X_n$ will forever jump back and forth along one singe edge (between $0$ and $1$, say). This is easily checked using the Borel-Cantelli lemma. Similarly, when $b=0$ and $a$ is negative, then with positive probability, the walk will be stuck on a set of two edges (say, it will visit only the three sites $0$, $1$ and $2$). These are the {\em purely self-attractive cases}.

\subsection {When $b>0$ and $-b/3 < a < b$: The TRSM regime?}

Let us first very briefly recall here the case where $b>0$ and $a=0$: The walker is then driven by the negative gradient of the local time at its current position.
As we have already mentioned, this is the so-called true self-repelling walk TSRW\footnote {In this terminology originating  in the physics literature, ``true'' refers to the fact that this is a true walk -- the distributions of the first $n$ steps is consistent with the distribution of the first $n+1$ steps -- and this is not always the case for measures on self-avoiding or self-repelling walks; the term ``myopic'' random walk has also been used.}. This scaling limit was clarified in \cite{T1} and \cite{TW}; it is the true self-repelling motion (TRSM).

\begin {figure}[htbp]
\begin {center}
\includegraphics [height=2in]{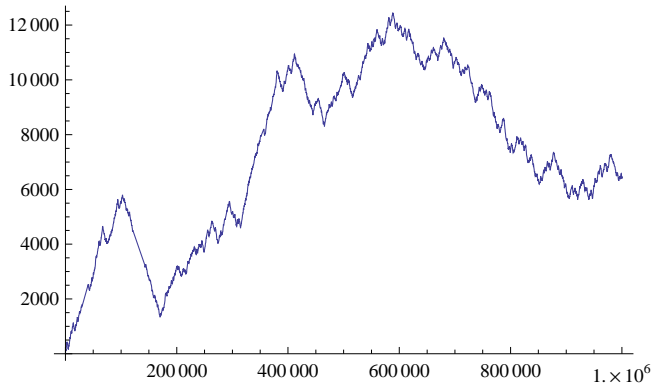}
\caption {TSRW trajectory: $a=0$, $b=1$, $P(x)= (x-1)$.}
\end {center}
\end {figure}

It is quite natural to expect that this TSRW scaling behaviour is stable under perturbation, i.e.  when choosing  $a\not=0$, but sufficiently small (in absolute value). There are some good reasons to guess that in  the parameter range $a \in (-b/3, b)$ (and for $b>0$), the long time asymptotic behaviour and scaling limit of the walk is similar to the TSRW case.
Indeed, besides numerical evidence, one can note that:
\begin {enumerate}
 \item
In this parameter range, we can identify a stationary and ergodic distribution of the ``environment seen by the random walker'' process.
This distribution is qualitatively similar to the TSRW case, the difference being that instead of a product measure we get an exponentially mixing Gibbs measure.
For details, see subsection \ref{Sstationary} below.
\item
In this stationary and ergodic regime, using the variational method initiated in \cite{lqsy} and exploited in \cite{ttv} for similar but continuous space-time models,  one
 could very likely obtain superdiffusive lower bounds of the type $E(X_n^2)\ge C n^{5/4}$ for the asymptotic variance of the displacement of the random walker.
For details of similar computations in continuous space-time, see \cite{ttv}.
\end {enumerate}

\subsection {When $b > 0$ and $a<-b/3$: The stuck case}

In this regime, the competition between the nearest edge self-repellence due to  $b > 0$ and next-to-nearest edge self-attractiveness due to  $a> 0$ is won by the latter because $|a|$ is large enough, and the random walker is eventually trapped in an interval of finite length. The length of the trapping interval increases to infinity as the ratio  $|a| / b$ decreases to the critical value $1/3$.
More precisely, let us define for all positive integer $k$,
$$
A_k:=1 + 2\cos\left( \frac{2\pi}{k+2} \right) .
$$
This is an increasing sequence ($A_1= 0$, $A_2= 1$, $A_4= 2$ and $\lim_{n \to \infty} A_n = 3$).
In \cite{ETW}, the following theorem is proved:
\begin{theorem}[\cite{ETW}]
\label{thm_sticky} Suppose that $a< 0 < b$.
\null
\begin {itemize}
 \item
If $b/|a| \in (A_{k}, A_{k+1})$ for some $k \ge 1$, then with positive probability, the walk remains stuck on a set of $k+2$ consecutive sites (and visits all these sites infinitely often).
\item
If $b/|a| > A_{k}$, then almost surely, the walk does not get stuck on a set with less than $k+2$ sites.
\end {itemize}
\end{theorem}

For instance, for some values of the parameter $|a|$ larger but very close to $b/3$, the walk can remain stuck on a set of exactly $10^9$ sites. On the other hand, when  $ 0 < -a  \le b/3$, then the walk does almost surely not get stuck on a finite set.

\begin {figure}[htbp]
\begin {center}
\includegraphics [height=1.7in]{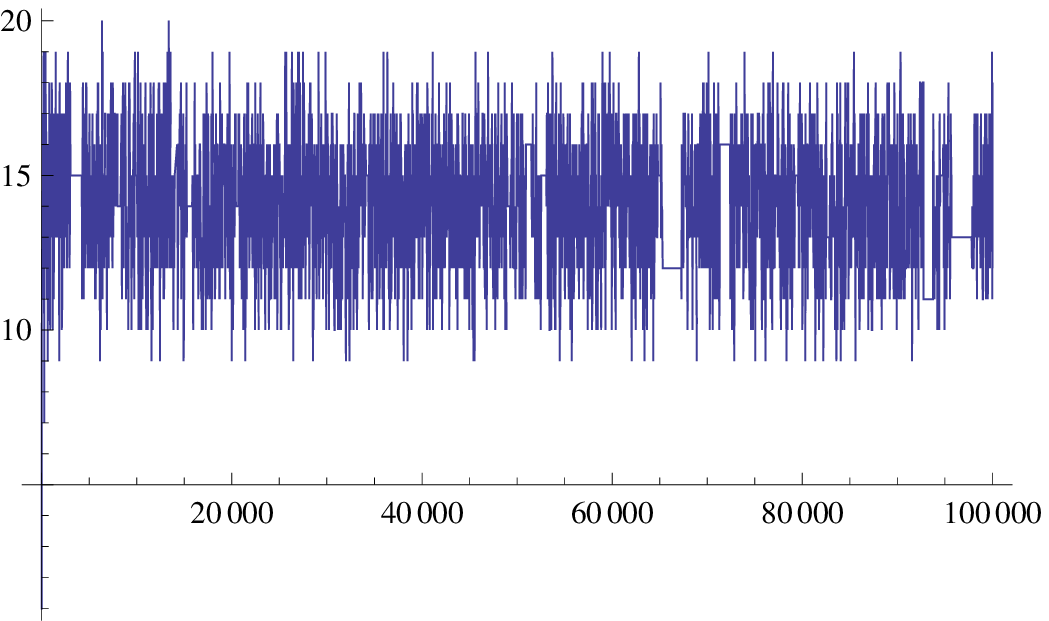}
\includegraphics [height=1.7in]{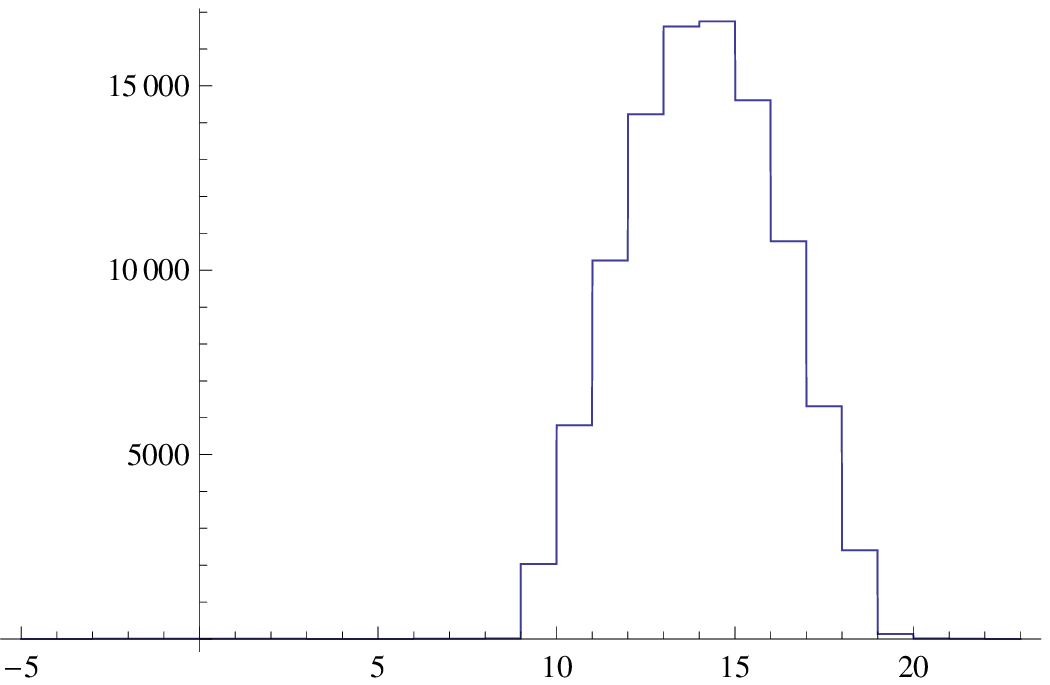}
\caption {A walk ($a=-1.1$, $b=3$) stuck on 11 sites: The trajectory and the local time}
\end {center}
\end {figure}
%We therefore see that the phase-transition between ``stuck'' and ``non-stuck'' occurs exactly when $a=-b/3$.
%One can note that this corresponds exactly to the case where $D$ is a discrete third derivative of $\ell$.
For proof of Theorem \ref{thm_sticky} and more refined details on the asymptotic behaviour of the walk in this particular trapping regime see \cite{ETW}.

\subsection{When $b \ge 0$ and $a>b$: The slow phase?}

In this regime, nearest and next-to-nearest edges are all repelling. Nevertheless these effects do not just always add up to produce a TSRM-type asymptotic behaviour. The second neighbours can produce traps that are difficult to get out from.

In order to illustrate this, let us describe  a scenario of building up traps, for the extreme case where $b=0$ and $a>0$:
The random walk starts by jumping back and forth on the edge $e_0$ between $0$ and $1$ a large even number $N$ of times. This happens with probability $2^{-N}=e^{-c_1 N}$. Then, in some $O(1)$ number of steps it moves to site $-4$.  (Note, that due to the large number of visits at edge $e_0$ the walker does not have a chance to return from site $-1$.) Now, it performs yet again a number of $N + o(N)$ back and forth jumps on the edge $e_{-4}$ (between sites $-4$ and $-3$). This happens with probability $e^{-c_2 (N+o(N))}$. After all this, in $O(1)$ number of steps the walker lands on edge $e_{-2}$ (between sites $-2$ an $-1$). Now, due to the  $N + o(N)$ number of visits at edges $e_0$ and $e_{-4}$ the walker does not have a chance to get away from edge $e_{-2}$ before jumping  $e^{aN + o(N)}$ times back and forth on this edge.

Hence, we see that while it is very unlikely to build up a trap (and it therefore typically takes a long time to build up one), once it is built and if $e^{a}$ is large enough, then it takes even much longer to escape from it. In fact, the previously described strategy to build such traps is certainly not the optimal one, and the probability to create one is in fact smaller, so that for all positive $a$, it takes much longer to escape from the trap than to create it. This slowing-down mechanism is reminiscent of that of random walks in random environment, and it indicates that there is no scaling limit process with continuous trajectories.

\begin {figure}[htbp]
\begin {center}
\includegraphics [height=1.7in]{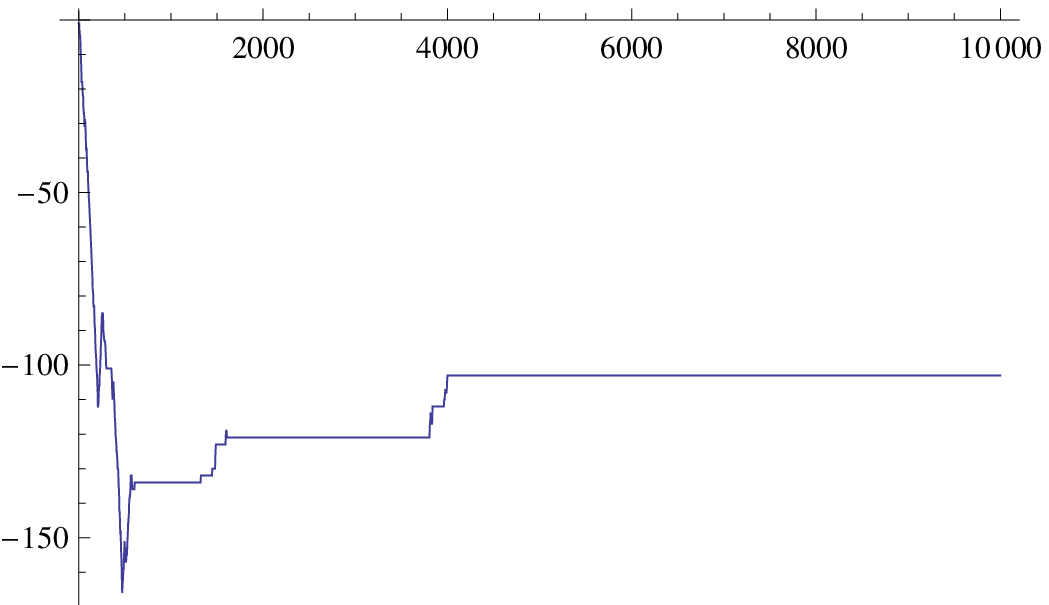}
\includegraphics [height=1.7in]{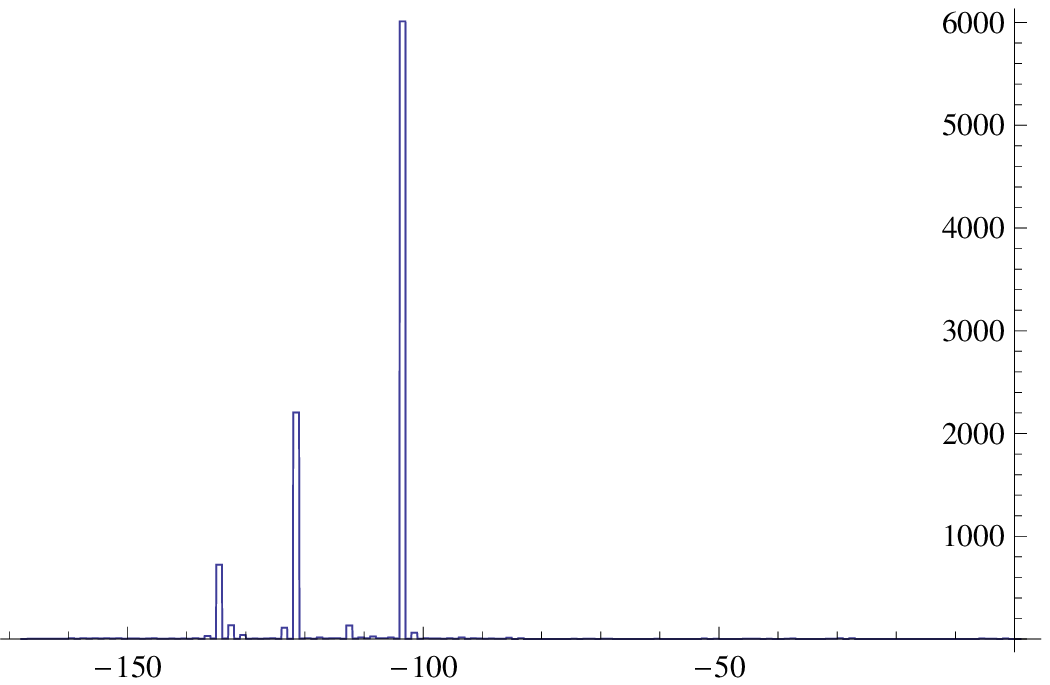}
\caption {The trapping phenomenon for $a=1$ and $b=0$: The trajectory and the local time}
\end {center}
\end {figure}

\subsection{The two critical cases}

\begin {figure}[htbp]
\begin {center}
\includegraphics [height=2in]{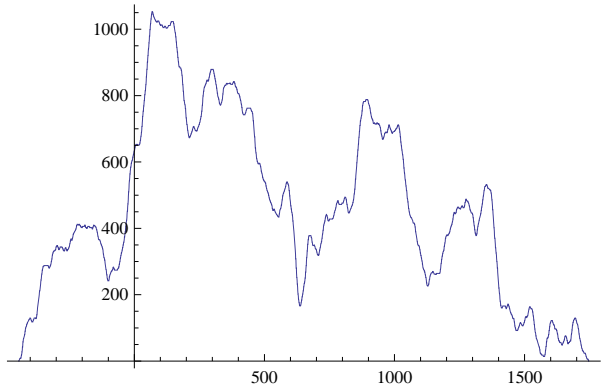}
\caption {The third derivative walk: The (smoother) local time profile}
\end {center}
\label {figlt3}
\end {figure}

When $b$ is positive, we do not have precise results at this point for the two borderline cases where $a=-b/3$ and $a=b$.
\begin {itemize}
\item $a= -b/3$  is the borderline case between the conjectural TSRM regime (where the walk should scale like $n^{2/3}$) and the sticky regime (when the walk is eventually stuck to a finite interval). We have seen that the walk does almost surely not stay in a finite interval, but its asymptotic behaviour is not really well understood. Naive scaling arguments (and the stationary measure computations mentioned at the very end of the next section) would at first suggest that the non-degenerate scaling of the walk should be like $n^{2/5}$. But it can be rigorously proved that in a stationary regime there exists a diffusive lower bound on the asymptotic variance of the displacement, and intuitively, one might guess diffusive behaviour (the stationary local time regime -- see also Figure \ref {figlt3} -- and the walk displacement might have different scaling exponents and not interact in the large-scale limit).
\begin {figure}[htbp]
\begin {center}
\includegraphics [height=2in]{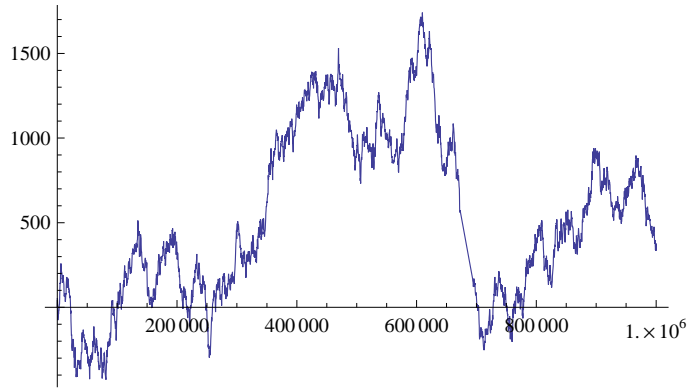}
\caption {The ``third derivative walk'' i.e. $P(x)=(1-x)^3$: The trajectory}
\end {center}
\end {figure}

\item The case where $a=b$ can be interpreted as the TSRW with self-repulsion defined in terms of occupation times on sites rather than edges. Its scaling behaviour is expected to be similar to the TSRW on edges (see
\cite {ETW} for references).
\end {itemize}

\begin {figure}[htbp]
\begin {center}
\includegraphics [height=3in]{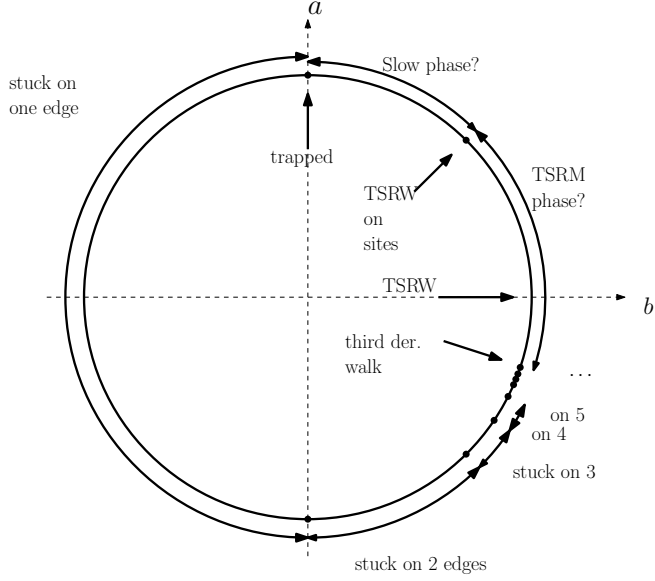}
\caption {The phase diagram in terms of $(b,a)$ on the unit circle}
\label {phased}
\end {center}
\end {figure}

Figure \ref {phased} sums up  the rough phase diagram of the asymptotic behaviour of the walk depending on the values of $(b,a)$ on the unit circle.

\subsection{Stationary measures for the cases where $b > 0$ and $-b/3<a<b$}
\label{Sstationary}

We now  quickly survey some computations related to the existence  of stationary measures for the current environment (i.e., the local time fluctuations) as seen from the random walk in the conjectural ``TSRM'' regime.
We will work in the more general context where the dynamics is driven by a height-invariant linear function of the local time i.e. where
$$ R ((\ell)) = e^{D} / ( e^D + e^{-D})$$
with
$$ D((\ell)) = \tilde D ((\eta)) = \sum_{x} \aba(x) \eta (x),$$
where all but finitely many $\aba (x)$'s are zero.
 We will write $ R ((\ell)) = p( \eta)$.

It is convenient to work with an initial profile $\eta_0$ that belongs to the set
$$ \Omega = \{\eta \ : \ \hbox  { for all } x, \ \eta(x) +1_{x=0} \in  2 \Z \}.$$
It is immediate to check that this parity condition will be conserved for all $n\ge0$ (the opposite parity convention i.e. $\eta +1 \in \Omega$ would also be preserved) so that
 $(\eta_n, n \ge 0)$ is a Markov chain on $\Omega$. We can also restrict ourselves to the family of $\eta$'s in $\Omega$ with finite
  norm $ \| \eta \| = \sum_{x \in \Z} 2^{-|x|} |\eta (x)| $.

We will use the following notation. If $X_{n+1} - X_n = 1$ i.e. $X$ jumps to the right, we write $ \eta_{n+1}= {\mathcal R} \eta_n$, and if it jumps to the left, we write
$\eta_{n+1}= {\mathcal L} \eta_n$.
In other words, for all $x$,
$$
{\mathcal R}\eta (x) = \eta (x+1) + \delta_0 (x) - \delta_{-1}(x),
\quad
{\mathcal L}\eta (x) = \eta (x-1) - \delta_0 (x) + \delta_{+1}(x).
$$
%The maps ${\mathcal R}^{-1}$ and ${\mathcal L}^{-1}$ are easily computed and are given by similar formulas.

The Markov chain is described by the fact that (given $\eta_n$) $\eta_{n+1}$ is equal to ${\mathcal R} \eta_n$ with conditional probability $p (\eta_n)$, and that $\eta_{n+1}= {\mathcal L} \eta_n$ otherwise. The transition operator of the Markov process is therefore
$$
Pf ( \eta) =  p(\eta) f({\mathcal R} \eta ) + (1- p(\eta)) f ( {\mathcal L}\eta).
$$

As often, it is convenient to view this Markov chain as the jump chain associated to a continuous-time Markov chain on the same state space. This continuous-time Markov chain $(\tilde \eta_t, t \ge 0)$ jumps to the right with intensity $dt \times \exp \{{\tilde D} (\tilde \eta_t) \} $
 and to the left with intensity
$dt \times \exp \{-{\tilde D} (\tilde \eta_t)\}$.

The particular choice for the dynamics leads us naturally to look for a stationary distribution for $\tilde \eta_t$ of ``discrete Gaussian type''. We are going to suppose in the present section that:
\begin {itemize}
 \item
$(\aba(\cdot))$ is even (i.e., $\aba (x)= \aba (-x)$ for all $x$). This corresponds to the left-right symmetry mentioned in the Introduction.
\item The self-interaction has finite spatial  range (i.e., all but finitely many $\aba (x)$'s are equal to zero)
\item
It is positive definite in the sense that for all $(\lambda (x), x \in \Z)$
such that all of them but a finite positive number is equal to zero,
\begin{equation}
\label {qform}
\sum_{x,y \in \Z} \lambda(x)\lambda(y)  \aba (y-x) >  0
\end{equation}
(this implies in particular that $\aba (0)>0$).
\end {itemize}
Note that in the case where $D=D^{a,b}$ and $b>0$, this last condition precisely means that $ -b/3 <  a  < b$.

\medbreak
Under these assumptions there exists a unique Gibbs measure $\pi_0$ on $\Omega$ formally given by
$$
d \pi_0(\eta)\sim
\exp \{ - \frac 1 2 \sum_{x,y} \eta(x) \eta (y) \aba(y-x) \}.
$$
The proper, rigorous definition of the Gibbs measure $\pi_0$ is in terms of conditional specifications on finite intervals, given fixed boundary conditions. This Gibbs measure is well understood, and  its correlations decay exponentially fast with distance and the zero-one law holds for tail events.

\begin{proposition}
Under the previous assumptions, the probability measure $\pi_0$ is stationary and ergodic for the Markov process $\tilde \eta_t$.
\end{proposition}

The proof is straightforward and based on the following observation:
$$
\frac{d\pi_0( {\mathcal R}\eta)}{d\pi_0(\eta)}
=
\exp\{-(\aba (1) - \aba(0)) - \sum_{x\in\Z}(\aba (x) - \aba (x-1)) \eta (x) \}
=
\frac {\exp \{ \tilde D (\eta) \} }{\exp \{\tilde D( {\mathcal R} \eta)\}}
$$
and similar formulae when one replaces ${\mathcal R}$ by ${\mathcal L}$. It proves that in fact the measure $\pi_0$ is invariant under the process with ``right-jumps'' only, and under the process with ``left-jumps'' only.

Ergodicity also follows from the above observation. The Dirichlet form
$$
{\mathcal D}(f)=\int_\Omega f(\eta)(I-P)f(\eta)d\pi_0(\eta)
$$
is easily computed. Using the above observation we get
$$
2{\mathcal D}(f)=
\int_\Omega \big(
p(\eta)(f(\eta)-f(\mathcal R\eta))^2
+
(1-p(\eta))(f(\eta)-f(\mathcal L\eta))^2\big) d\pi_0(\eta).
$$
Hence, since $0<p(\eta)<1$ almost surely, it follows that $\mathcal D(f)=0$ if and only if almost surely
$$
f(\mathcal R\eta)=f(\eta)=f(\mathcal L\eta).
$$
But it is standard fact that the above almost sure identities hold if and only if $f$ is tail measurable. Hence, by the zero-one law mentioned before, $f$ is almost surely constant. This proves ergodicity.

\medbreak

Also, one can note that the mean waiting time of the continuous-time Markov chain at $\eta$
is given by 
$$
\tau (\eta) := e^{\tilde D(\eta)}+e^{-\tilde D(\eta)}, 
$$
and that $Z := \pi_0 (1/ \tau)$ is finite. From this, it follows that one can define a probability measure $\pi$ on $\Omega$ by
$$
\frac {d \pi}{d \pi_0} (\eta) = \frac 1 {Z \times \tau (\eta)}
$$
and that:

\begin {corollary}
 This probability measure $\pi$  is stationary for the discrete chain $(\eta_n)$.
\end {corollary}

Hence, we see that in this regime, and under these well-chosen starting distributions, the large-scale behaviour of the local-time profile is of Brownian-type.

We can note that in the borderline case where $b>0$ and $a= -b/3$ (i.e. the discrete ``third derivative'' case), the previous stationary distribution analysis can be adapted to show that the discrete gradient of $\tilde \eta_t$ (that is: the discrete second derivative of the local time profile) has a stationary measure of discrete Gaussian type.
More precisely, we see in this interesting case that if we define for all $e \in E$,
$$ \nabla \eta (e) = \eta (e+1/2) - \eta (e-1/2), $$
then the process $( \nabla \tilde \eta_t, t \ge 0)$ is a Markov process, and that it has a stationary distribution given by i.i.d. discrete Gaussians.

For details of similar computations of stationary and ergodic measures of the environment seen by the moving particle, in continuous space-time setup, and their consequences see also \cite{ttv}.

\subsection {Some comments}

To conclude this section, let us very quickly list a few observations:
\begin {itemize}

\item
In the setting of the TSRW, it was natural to describe the dynamics of the walker by saying that it was driven by an approximate negative gradient of the
local time at its current position, and the dynamics of its scaling limit could also be described in this way. Our previous analysis however shows that this
type of description has to be handled with care. The case $\{ a=1, b=0 \}$ seems at first to also correspond
to a negative discrete gradient, but its large scale behaviour is completely
different.

\item
In fact, in the general case, what seems more relevant to classify the asymptotic behaviour of the walk is the quadratic form appearing in (\ref {qform}) or alternatively the study of the roots of
$P (x)$. This is also apparent in the study of the stuck walks in \cite {ETW}. Note also that the borderline cases correspond to (multiples of) the polynomials $(1-x)^3$ and $(1-x)^2(1+x)$.  In the general case (with larger interaction-range), the phase diagram is seemingly more complicated, with different connected components  of the parameter space corresponding to similar asymptotic behaviours.

\item
A particular role is then played by the ``odd'' derivatives, i.e. more precisely, the polynomials $(1-x)^{2n+1}$. Then, the discrete $n$-th derivative of $\tilde \eta$ (and therefore of $\eta$) has a stationary measure of discrete Gaussian type.
\end {itemize}

\section {Some cases without left-right-symmetry}
\label {S3}

Let us stress that (as opposed to the cases studied in the previous section), the examples studied in the present section break the left-right symmetry i.e. the law of $(X_n, n \ge 0)$ and $(-X_n, n \ge 0)$ are different even if one chooses $L_0$ to be symmetric.

\subsection {Setup and statement}

In this section, we will first focus on the case where $D = - \Delta /2 $, where
$$ \Delta ( (\ell(e), e \in E)) = \ell(-3/2) - \ell(-1/2) - \ell(1/2) + \ell(3/2).$$
 We then choose
$$ R (( \ell (e), e \in E)) = \frac { e^{-\Delta/2}}{ e^{\Delta/2} + e^{-\Delta/2}},$$
which penalizes jumps in the direction of this second derivative. The precise form of the function $R$ is not really important for what follows. For instance, our results clearly remain true if we only assume in addition of locality and height-invariance, that
\begin {itemize}
 \item $R ((-\ell)) = 1 - R ((\ell))$.
 \item $R((\ell))$ is a function of $\Delta$  goes to $0$ very quickly as $\Delta \to +  \infty$.
\end {itemize}
These conditions can in fact also be relaxed.

Note that $\Delta= \eta (1) - \eta (-1)$
 is a discrete and symmetric version of the second derivative of the local time, but this feature is not essential in what follows.
Intuitively speaking, the walk will tend to move to the right if the local time profile is locally concave and to the left if it is locally convex.

The main result of the present section is the following: Suppose that the initial profile is identically zero (again, this is not really necessary -- we only need this initial profile not to be too ``wild'').
Then:

\begin {proposition}
\label {t2}
With positive probability, $X_n \sim \sqrt {2n}$
 as $n \to \infty$.
\end {proposition}

On large scale and on a positive fraction of the probability space, the walk becomes almost deterministic, but non-ballistic, and the scaling is of Brownian type.
 As we shall see, the proof gives in fact a more precise description of the behaviour of $(X_n, n \ge 0)$ on this event of positive probability.

\begin {figure}[htbp]
\begin {center}
\includegraphics [height=1.7in]{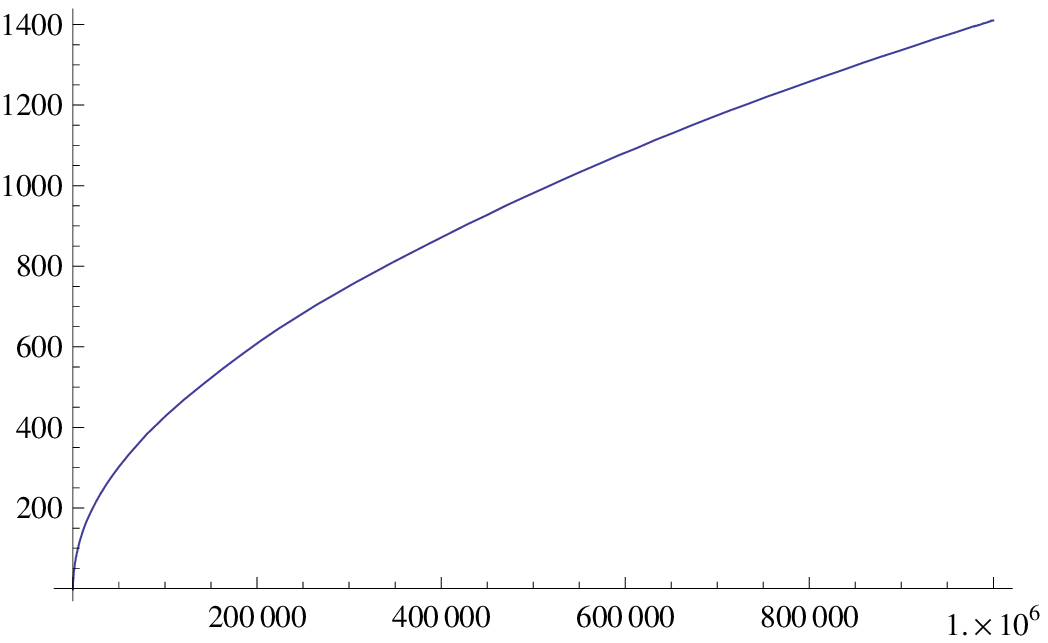}
\includegraphics [height=1.7in]{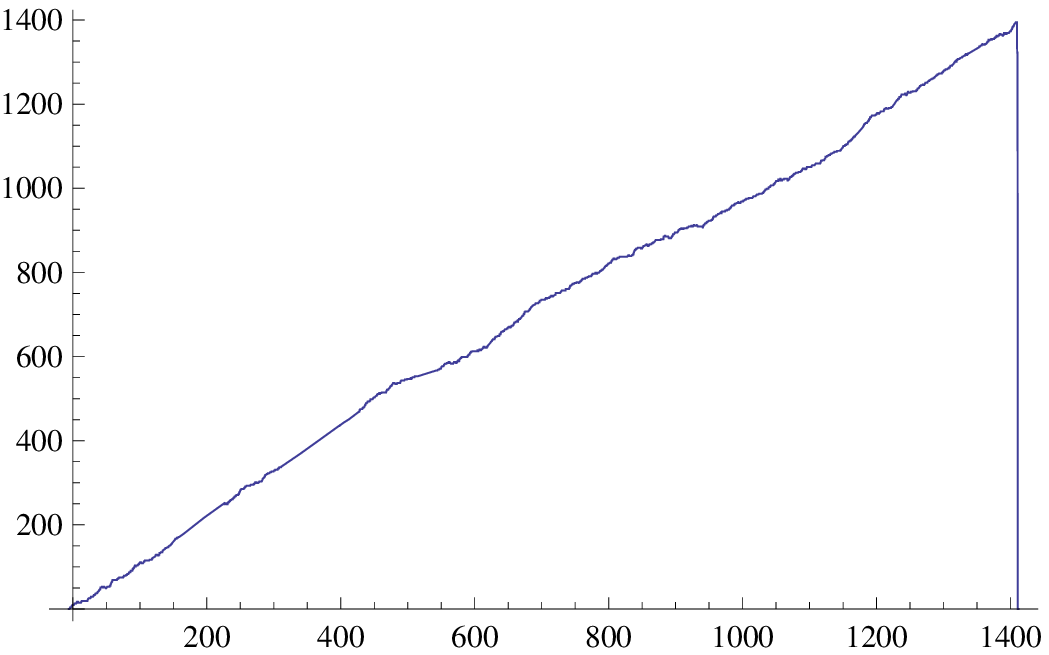}
\caption {The second derivative walk: The trajectory and the local time}
\end {center}
\end {figure}

Let us immediately stress that the proof will in fact show that the results remain true if one replaces $\Delta$ by any function $\tilde \Delta$ of the type
$$ \tilde \Delta ( (\ell (e), e \in E)) = \ell(-3/2) - \ell(-1/2) - \ell(1/2) + \sum_{k \in E, \ k \ge 3/2} a_k \ell(k).$$
The facts that $\Delta$ is a local version of the second derivative, that it is translation-invariant and that it is ``anti-symmetric'' are therefore not so crucial after all.
In fact, our conclusion holds true if we replace $\Delta$ by any local function
$\hat \Delta ( (\ell (e), e \in E)) = \sum_{ e \in E} a_e \ell (e) $
such that
\begin {itemize}
\item
For some $e_0 < 0$, $a_e =0$ for all $e \le e_0$.
\item
$\sum_{e < 0, e \in E} a_e = 0 $
\item
$a_{-1/2}$ and $a_{1/2}$ are both negative.
\end {itemize}
This for instance includes discrete $p$-th derivatives for ``even'' $p$'s such as
$$\hat \Delta ((\ell))
= - \ell (-5/2) + 3 \ell (-3/2) - 2 \ell (-1/2) - 2 \ell (1/2) +  3\ell (3/2) - \ell (5/2)$$
that corresponds to a discrete forth derivative.

In the ``continuous limit'', the occupation time density profile of the motion $t \mapsto \sqrt {2t}$ has a discontinuity at its local position (i.e. at time $t$, it is equal to $x 1_{x \in [0,\sqrt 2{t}]}$), and the size of  this jump is growing with time. Hence, one could a priori naively expect that the walk would not slow down, but on the contrary speed up when this gap widens. This can be explained by the fact that the self-interacting random walk therefore captures some feature of the discrete local time profile that is not visible in the continuous limit.

\subsection {The scenario}
The coming three subsections are devoted to the proof of Proposition \ref {t2}.
Let us first describe in plain words a possible behaviour for our walk $X$ (we shall then prove that indeed, this behaviour occurs with a positive probability):
The walk $X$ starts at the origin, jumps to its right i.e. $X_1 =1$, and it forever remains positive i.e. $X_n > 0$ for all positive times. In fact, when $n \to \infty$,  $X_n$
goes off to $+ \infty$. For each $n \ge 0$, denote
$S_n = \max ( X_m, m \le n)$ the past maximum of the walk.
In our scenario, $X_n \ge  S_n -  2$ for all $n$.
In other words, for each $n$, either $X_n$ is equal to its past maximum $S_n$, or it is equal to $S_n -1$ or $S_n -2$. 

We see that we can therefore decompose the set of times as follows: Define for each $x \ge 0$ the hitting time $\tau_x$ of $x$ by $X$, and let $I_x = [\tau_{x+1}, \tau_{x+2})$.
Then, in our scenario, during each of the intervals $I_x$, the walk jumps back and forth on the edges between $x-1$ and $x+1$. When it stops doing so, it first jumps to $x+2$, and then jumps back and forth on the edges between $x$ and $x+2$ and so on.

Imagine for a moment that this has happened for a while and that just before $\tau_x$ (for a large $x$) i.e. when $n = \tau_x -1$, the following loosely defined event $G(x)$ is true:
The walker has visited each of the edges $\{x-2, x-1\}$,
 $ \{x-3, x-2\}$ and $\{x-4, x-3\}$ many times, and the number of times it has visited these three edges remain however comparable (in the sense that the difference between these three numbers of visits is small). In particular, immediately before it chooses to jump to $x$ for the first time, the value of $\Delta((\ell_n)) 
 = \ell_n (-3/2) - \ell_n (-1/2)$ is rather close to zero (mind that at that moment $X_n=x-1$, and  that $\ell_n (1/2)= \ell_n (3/2)=0$ because the walk has not visited $x$ yet) so that the probability to indeed jump from $x-1$ to $x$ is therefore not too small.

  Then, the walk arrives at $x$ for the first time. At that moment, the two edges to its right have not yet been visited, the edge between $x$ and $x-1$  (to its left) has been visited only once (because the walker arrives at $x$ for the first time). On the other hand, the edge between $x-1$ and $x-2$ has been visited a lot of times. Hence, the walk will (very likely) jump back to $x-1$. Once it is back at $x-1$, the probability to jump to $x$ or to $x-2$ at that moment is neither very small nor very close to $1$ (because the situation can not drastically change because of the two previous jumps). Note that if the walk then jumps to $x-2$, the value of $\Delta$ will then be very small, so that the walk will want to jump immediately back to $x-1$. Hence, the walk is (with high probability) trapped between $x-2$ and $x$ for a while. How long does this happen? Well, one should note that each time the walk jumps on $\{x-2, x-1\}$ or on $\{x-1, x\}$, it will increase the chance to jump to the right next time it is at $x-1$. Hence, after a short while, the walk will in fact only jump back and forth between $x-1$ and $x$. This will be the case until the number of times at which it has jumped on $\{x-1, x\}$ starts to be comparable with the number of times at which it has jumped on $\{ x-2, x-1\}$, because the walk will then have a significant chance to move to $x+2$. Hence, we end up in a situation where $G(x+1)$ holds.

\begin {figure}[htbp]
\begin {center}
\includegraphics [width=2in]{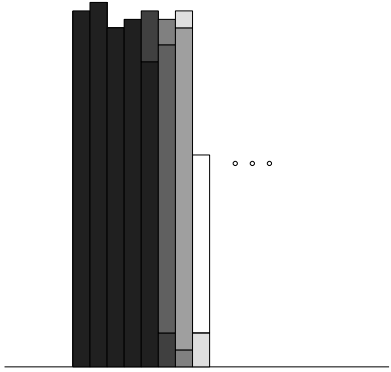}
\caption {The local time profile in our scenario: Darker cells are constructed first.}
\end {center}
\end {figure}

 One important point in this scenario is that in the end, the walk will tend to jump a little more often on $\{ x, x+1\}$ than on $\{x-1, x \}$. So, by a law of large numbers type argument,   the total number of visits of the edge $\{x-1, x \}$ will turn out to grow ``linearly'' in $x$ as $x \to \infty$. This will a posteriori justify our assumption (in the definition of $G(x)$) that at the time $\tau_x$, the walk has visited $\{x-2, x-1\}$ and $\{x-3, x-2\}$ many times.

\subsection {Auxiliary sequences}

We now start to turn the previous ideas into a rigorous proof.
Let us define a family of independent random variables $(\xi_{j,k}^{x})$, where $j \in \Z/2$, $k \ge 1$, $x \in \Z$. For each $j$, $k$ and $x$, the law of the random variable is the following:
$$
P ( \xi_{j,k}^x = +1 ) = 1 - P ( \xi_{j,k}^x = -1 ) = \frac {e^{-j}}{e^{j} + e^{-j}}.
$$
We can then define (deterministically) our random walk $X$ using this family of random variables as follows.
Let $X_0 = 0$ and for each $n \ge 0$, choose inductively
 $$X_{n+1} = X_n + \xi_{\delta,k}^x,$$
where $n$ is the $k$-th time  at which $X_n=x$ and $\Delta_n=2\delta$. Note that when $X_n \ge 2$, then $\Delta_n$ is necessarily even (because $\ell_n (-3/2)$ and $\ell_n (-1/2)$ are odd, while 
$\ell_n (1/2$ and $\ell_n (3/2)$ are even). 
We will for the time being forget about this coupling between $X$ and the $\xi$'s, and just list
 a few simple observations concerning this
family $(\xi_{j,k}^x)$ of random variables.

A particular role will be played by the random variables $\xi_{j,k}^x$ when $x \ge 0$, $k=1$ and $j$ is an integer. We denote them by $\xi_j^x$ (and drop the subscript $k$). Let us insist on the fact that $\xi_j^x$ are not defined for semi-integer $j$'s (as opposed to the $\xi_{j,k}^x$'s). Let us also keep in mind that
for a given $k$ and $x$, when $j$ is very large, the probability that
$\xi_{j,k}^x = -1$ and $\xi_{-j, k}^x = +1$ is very close to $1$.

We now define certain events:
\begin {itemize}
\item
Let $A$ denote the event that for all $x$, all $j > \sqrt {|x|}/100$ and all $k \le 100 x^2$, $$\xi_{j,k}^x = - \xi_{-j,k}^x = -1.$$
\item
Note that almost surely, for each $x$, the number of positive (respectively, negative) $j$'s for which
$\xi_j^x$ is positive (respectively, negative) is finite (it follows immediately from the definition that it has a finite expectation).
 $B$ is the event that for all $x \ge 0$,
$$ \sum_{j  > 0} \left( 1_{\xi_j^x = 1 } +  1_{\xi_{-j}^x = - 1 } \right) \le \sqrt {x}/100.$$
\item
For all $x \ge 0$, let us define
$$
M^x = \sum_{j \ge 0} 1_{\xi_j^x = 1 } - \sum_{j < 0} 1_{\xi_j^x = - 1 }.
$$
There is no problem with the definition of $M^x$ for the same reason as above.
Note that because of symmetry, $ E( M^x ) = P ( \xi_0^x = 1 ) = 1/2$. The law of large numbers therefore  ensures that almost surely
$$ \lim_{y \to \infty} y^{-1} \sum_{x=0}^y M^x  =  1/2.$$
We define the event $C$, that for all $y \ge 1$, $ \sum_{x=0}^y M^x \in [y/4, 4y ]$.
\end {itemize}

A simple application of the Borel-Cantelli Lemma enables to prove that:

\begin {lemma}
With positive probability, $A \cap B \cap C  $ holds.
\end {lemma}

The proof is elementary and safely left to the reader.

\subsection {The coupling}

In order to simplify our notations, we will assume that at time zero, the initial profile is such that
$ \ell_0 (-1/2)= \ell_0 (-3/2) =1$ and that $\ell_0$ is zero otherwise. Then, in our ``good scenario'',  the $\Delta_n$'s will always remain even (because $\ell_n (-3/2)$ and $\ell_n (-1/2)$ are both odd, whereas $\ell_n(1/2)$ and $\ell_n (3/2)$ are both even).

This assumption is not a problem. If we start with $\ell_0$ being identically equal to $0$, we let $X$ move twice to the right, and at this time $n=2$, the situation is exactly the previous one. Since we will prove asymptotic results about $X$ with positive probability, our theorem will follow immediately from the results with this particular initial profile.

Recall that our walk $(X_n, n \ge 0)$ is defined deterministically from the $\xi_{j,k}^x$'s as follows: $X_0=0$ and for each $n \ge 0$,
$$ X_{n+1} = X_n + \xi_{j,k}^x $$
where $x = X_n$, $j = \Delta_n / 2$ and $k$ is the cardinal of the set
$$\{m \le n \ : \  X_m = x \hbox { and } \Delta_m = 2j \}.$$
Clearly, this definition ensures that one uses never the same $\xi_{j,k}^x$ twice, and that indeed
$$ P ( X_{n+1} = X_n + 1 \ | \ X_0, \ldots, X_n ) =
\frac {e^{-\Delta_n/2}}{e^{\Delta_n /2}  + e^{- \Delta_n /2 }}
.$$

We are going to assume that  $A \cap B \cap C $ holds, and we will show that our ``good scenario'' holds as well. Note that there is nothing probabilistic in the following arguments.
For this, we will inductively prove that for each $x \ge 0$, certain events $E(x)$ hold. Define for each
$x \ge 1$, the stopping time $\sigma_x= \sigma (x)$ as the first time at which the edge $\{x-1, x\}$ has been crossed more than $x/8$ times and the walk is at $x-1$.
We say that the event $E(x)$ holds if the following four conditions hold:
\begin {enumerate}
\item
The walk did not visit $x+1$ before $\sigma_x$.
\item
On $[\sigma_{x-1}, \sigma_x ]$, the walk did only visit the three sites $x-2, x-1$ and $x$.
\item
The last $x/10$ jumps of the walk before $\sigma_x$ were all on the edge between $x-1$ and $x$.
\item
At the time $n=\sigma_{x}$, the two quantities $\ell_n (-1/2) $ and $\ell_n  (-3/2) $ differ by not more than $\sqrt {x}$, and they are both larger than $x/6$ and smaller than $50x$.
\end {enumerate}

\begin {figure}[htbp]
\begin {center}
\includegraphics [width=3in]{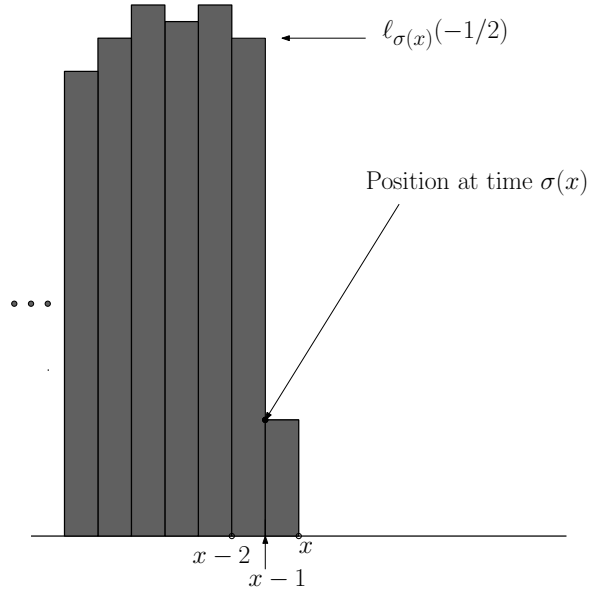}
\caption {The local time profile at $\sigma(x)$}
\end {center}
\end {figure}

Assume that $X_1=1$, $X_2=2$,  and that $E(3), \ldots, E(x)$ do hold (and that $A \cap B \cap C$ hold as well).

For simplicity, let us first assume that $\xi_{j}^x$ are all equal to $+1$ when $j \le 0$ and to $-1$ when $j > 0$ (we will then later see what difference it makes if for finitely many values of $j$ this is not the case). Because $E(x)$ holds, it implies that at the ca. $x/10$ first times at which the walk has been
at $x$ before time $\sigma_x$, the value of $\Delta_n$ at those times was greater than $x/20$ (because  $\ell_n(-1/2)$ was small, $\ell_n(-3/2)$ was large and
$\ell_n (1/2)= \ell_n (3/2) =0$).
Furthermore, the value of $\Delta_n$ at the time $n=\sigma(x)= \sigma$ is very negative (recall that $X_{\sigma (x)}= x-1$, that at this time $\ell_\sigma (-3/2)$ and $\ell_\sigma (-1/2)$ are comparable, while $\ell_\sigma (1/2)$ is greater than $x/8$ and $\ell_n (3/2)=0$).
Hence, our assumption on the $\xi_{j,k}^x$'s ensures that $X_{\sigma+1} = x$, $X_{\sigma+2}   =x-1$ and $X_{\sigma+3}=x$ (note that the $k$ used in those steps can not exceed the number of times one visited the sites). After these two steps, the value of $\Delta_n$ (when the walk is at $x-1$) has decreased by $2$ i.e.
$\Delta_{\sigma+2} = \Delta_\sigma -2$, whereas $\Delta_{\sigma+3} = \Delta_{\sigma+1} - 2 $ if $X_{\sigma+2} = x$.
Hence, we see that the walk will have to jump back and forth a number of times on the edge between $x-1$ and $x$; this will stop being the case because at some point, the walk will be at $x$ and will meet a $\xi_{j,k}^x$ that is equal to one. Because before that, the $\Delta$'s (at the visiting times of $x$) have been decreasing, we see that this can only happen when one uses a $k=1$ i.e. when one uses  $\xi_0^x=1$ (this is the $\xi_j^x$ that is equal to one, for which $j$ is the largest).

From that moment (call it $\tau$) onwards, the walk will start jumping back and forth on the edge between $x$ and $x+1$ (this is because when it is at $x$, it will use $\xi_j^x$'s for negative $j$'s, and when it is at $x+1$, it will use $\xi_j^{x+1}$ for very large values of $j$). This will happen until the time $\sigma_{x+1}$.
At that time,  $\ell_{\sigma(x+1)} (-1/2) = \ell_{\sigma(x)} (-1/2)$ i.e. the number of times at which the walk jumped on $\{x-1, x\}$ is equal to the number of times the walk jumped on $\{x-2, x-1 \}$, which ensures that in this very particular situation, $E(x+1)$ indeed holds.

What is the difference due to those $\xi_j^x$'s that are equal to $-1$ if $j \ge 0$ or to $+1$ if $j \le 0$? The first remark is that all these particular $j$'s will indeed be used in our process (note that $\Delta$ is decreasing two by two at each visit of $x$, and that it starts from a very positive value and become very negative, and therefore has to use all $\xi_j^x$'s for $j \in [ - \sqrt {x}, \sqrt {x}]$. Each time it meets a positive $j$ such that $\xi_j^x = +1$, it jumps to the right instead of jumping to the left. At the end of the day (i.e. at $\sigma(x+1)$) this will mean that $\ell(-1/2)$ will be diminished by two. Conversely, if it meets a non-negative $j$ such that $\xi_j^x = -1$, this will add $2$ to the number of jumps on $\{x-1,x\}$.
Hence, we see that in our ``real'' case,
$$ \ell_{\sigma (x+1)} (-1/2) = \ell_{\sigma(x)} (-1/2) + 2 M^x$$
and that  $E(x+1)$ will still hold because of our assumptions on the sum of the $M^x$'s.

As a consequence, we see that on the event $A \cap B \cap C$ of positive probability, all $E(x)$'s hold. This implies in particular that for each $x$, the walk will not come back to 
the site $x-1$ after $\sigma (x+1)$ i.e. that $\ell_{\sigma (x+1)} (-1/2)$ is in fact the total number of times the walk does jump on the edge
$\{x-1, x\}$.

Recall that $E(M^x) = 1/2$ and that almost surely, $\sum_{x=0}^y (2 M^x) \sim y $ as $y \to \infty$, so that on our event of positive probability,
$\sum_{x=0}^y \ell_{\sigma (x)} (-1/2) \sim y^2 /2$.
It follows that on our event of positive probability,
$ \sigma (x) \sim x^2/2$ so that $n \sim (X_n)^2 /2$  (recall that between $\sigma (x)$ and $\sigma (x+1)$ we know that $X_n$ is equal to $x-1$, $x$ or $x+1$) and $X_n \sim  \sqrt {2n}$ as $n \to \infty$.

This concludes the proof of Proposition \ref {t2}.

\subsection {An example with logarithmic behaviour}

We now focus on an example that exhibits yet another possible asymptotic behaviour for the walk: $D = - \Delta^\#$ where
$$\Delta^\# ( (\ell(e), e \in E)) = 2\ell(-3/2) - \ell(-1/2) - \ell(1/2) = -2 \eta (-1) - \eta (0)$$
and we then choose $F$ as before.  A naive first
guess would be that this walk is self-attractive (it is ``driven'' by the positive gradient of its local time) and that it should get stuck.
However:

\begin {proposition}
\label {p4}
With positive probability, $X_n \sim (\log n)/(\log 2)$ as $n \to \infty$.
\end {proposition}

\begin {figure}[htbp]
\begin {center}
\includegraphics [height=1.8in]{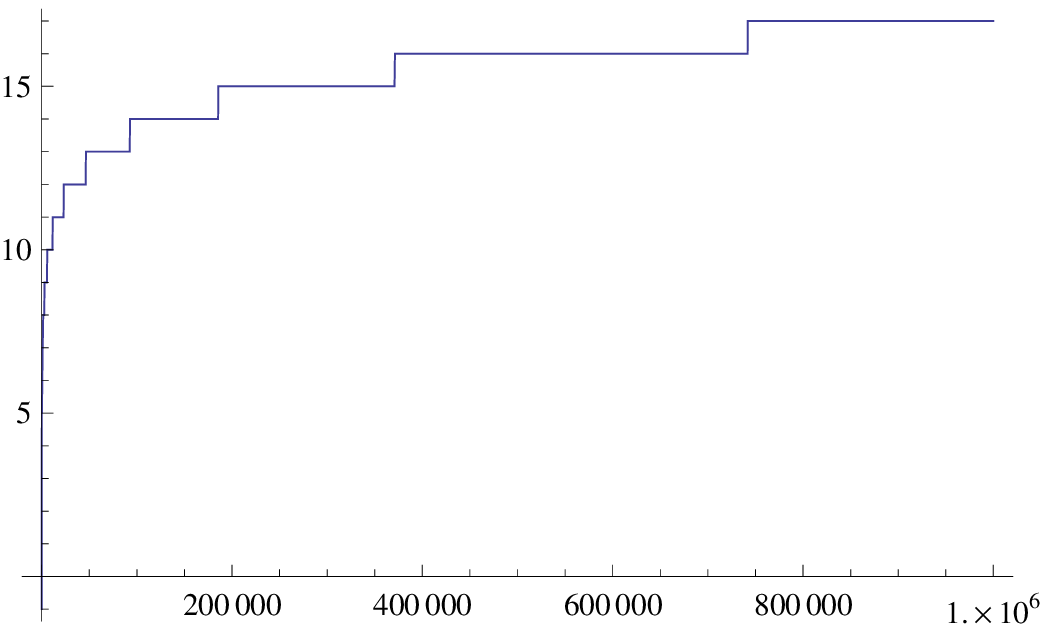}
\includegraphics [height=1.8in]{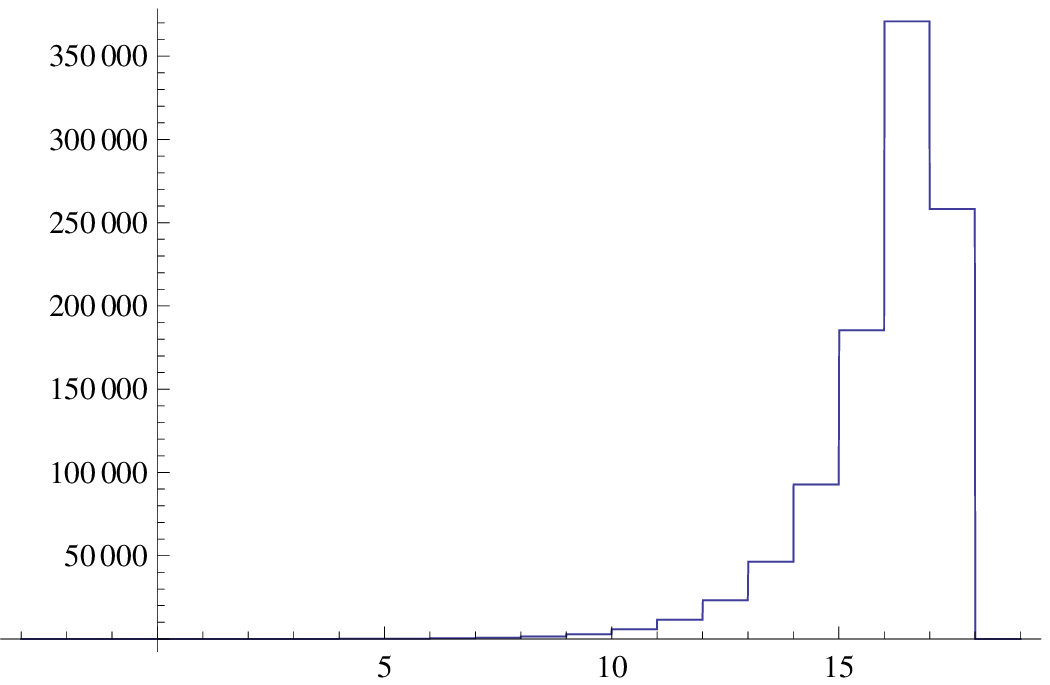}
\caption {The logarithmic walk: The trajectory and the local time}
\end {center}
\end {figure}

Again, the proof shows that this behaviour is valid for a wider class of self-interactions. In fact, it is quite similar to the case (with $\Delta$) studied in the previous subsections.
The main difference is that (in the ``good'' scenario that we will describe) the walk will visit approximatively twice more the edge $e+1$ than the edge $e$ for large $e$. Since the proof is otherwise almost identical to the previous case, we will only describe in plain words this ``scenario'', and leave the details of the proof to the interested reader.

Suppose that for some particular large time $n_x$:
\begin {itemize}
\item The walk is at its past maximum: $X_{n_x} = \max_{m \le n_x} X_m$. We call this site $x=X_{n_k}$, and for all $n \ge 0$, we define $U_n$, $V_n$ and $W_n$ to be
the respective values of $L_{n_x+ n}$ at $x-5/2$, $x-3/2$ and $x-1/2$.
\item The values of $U_0$, $V_0$ and $W_0$ satisfy:
$$ U_0 \in [(9/10) \times V_0/2,  (11/10) \times  V_0 / 2] \hbox { and } W_0 \in [V_0 /2 , V_0 ].$$
\end {itemize}
We also suppose that the value of $V_0$ is very large.
Note that under these assumptions, it very likely that $X_{n_x +1 } = x-1$ and $X_{n_x+2} = x$, i.e. that the walk will jump back and forth along the edge $x-1/2$ between $x$ and $x-1$,
because $ -2 U_0 + V_0 + W_0 $ is very large, while $-2V_0 + W_0$ is negative and has a large absolute value.
A quick analysis shows that (with high probability) the walk will jump back and forth on this edge until
the negative drift at $x$ stops being huge, and the walk will then jump to $x+1$ for the first time. When this happens, this means that at this time,
$2V_n$ and $W_n$ are comparable.
Then, for some time, it will jump on the three sites $x-1$, $x$ and $x+1$, but while doing so, the drift at $x$ (i.e., that it feels when it is at $x$) grows fast, so that
it will quickly be forced to jump along the edge between $x$ and $x+1$ only. We then let $n_{x+1}$ the first time $n$ at which
$ L_{n} ( x+1/2 )$ is greater than  $ L_{n} (x-1/2) / {2}$
and note that $n_{x+1}$ satisfies the same conditions as thos we required for $n_x$.
Note also that  once $n_x$ is found, the scenario  is very likely to hold until $n_{x+1}$ when $V_0$ is large (its conditional probability goes to $1$ very quickly with as $V_0$ gets larger). From this, it follows easily that with positive probability, the ``good scenario'' will indefinitely repeated. Then, clearly the total number of visits to the edge $y+1/2$ grows like $2^{y+ o(y)}$ as $y \to \infty$, and the proposition follows.

\begin {figure}[htbp]
\begin {center}
\includegraphics [width=2.5in]{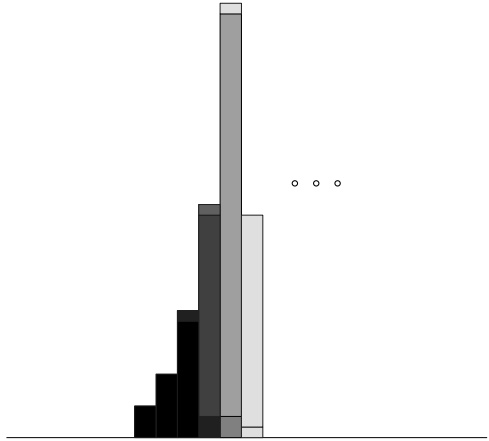}
\caption {The local time profile in the good scenario: Darker cells are constructed first.}
\end {center}
\end {figure}

\subsection {Ballistic behaviour}

To conclude, let us finally mention without proof
 a much less surprising possible asymptotic behaviour, namely the ballistic one. The following pictures correspond to the trajectory and local time for the case where $D = - \Delta^b$ where
$$\Delta^b ( (\ell(e), e \in E)) =  - \ell (3/2) -  \ell(-1/2) + 2 \ell(1/2)  = 2 \eta (1) + \eta (0)$$
that could seem at first glance to be again one of the walks driven (like TSRW) by the negative gradient of $\eta$. This example, like many others in the present paper illustrates how sensitive the discrete model is to little shifts in the definition of the driving dynamics.

\begin {figure}[htbp]
\begin {center}
\includegraphics [height=1.7in]{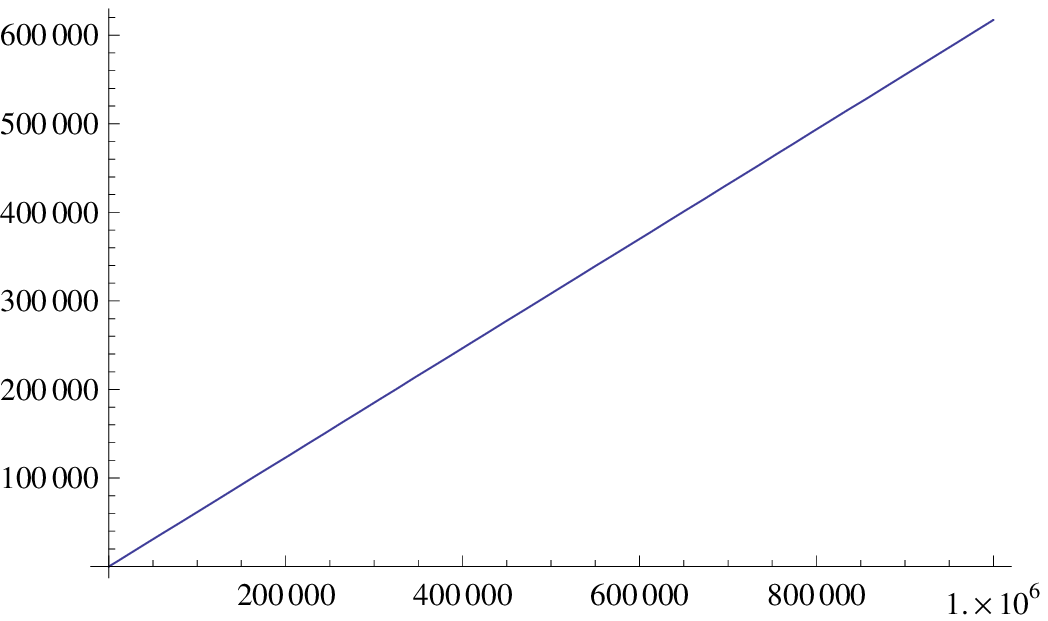}
\includegraphics [height=1.7in]{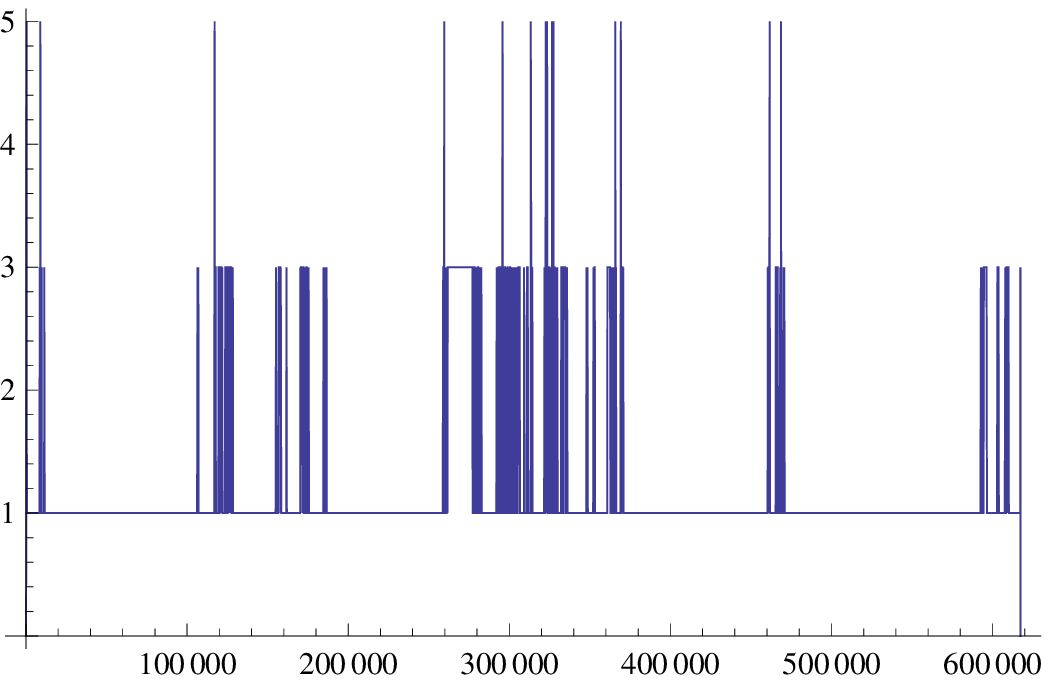}
\caption {A ballistic walk: The trajectory and the local time}
\end {center}
\end {figure}

\section {Some open questions}

Let us conclude with a list of (possibly accessible) open problems related to the questions that we have just discussed:

\begin {enumerate}

\item 
Prove the $2/3$ scaling behaviour of the walk in the ``TSRM regime'' (i.e. where $b>0$ and $a \in (-b/3,b)$ and the stationary distribution makes sense),
or even the convergence to TSRM. Or construct anote related Markovian model (for instance with another function $R$) where one can prove this?

\item 
Describe in some detail the ``actual'' (i.e., the one that actually dominates) trapping strategy in the ``trapped case'' $b=0$ and $a=1$.

\item 
Get some (even partial) description of the dynamics in the ``slow phase.'' Does the qualitative behaviour depend only on $b/a$?

\item 
In the case where the walk goes deterministically to infinity, are other scaling behaviours than $ct$, $c \sqrt {t}$ and $c \log t$ possible?

\item 
Improve some of the results to almost sure statements (instead of ``with positive probability'').

\item 
Is it possible that for some choice of the parameters in such  self-interacting random walks with finite range, the qualitative asymptotic behaviour is actually not almost sure (i.e. can two different behaviours that exist with positive probability exist?).

\item
Do qualitatively really new asymptotic behaviours arise when one considers larger (but finite) self-interaction ranges?

\end {enumerate}
Other questions more directly related to the ``stuck case'' are listed in \cite {ETW}.

\bigbreak

\noindent{\bf Acknowledgements.}
BT thanks the kind hospitality of Ecole Normale Sup\'erieure, Paris, where part of this work was done. The research of BT is partially supported by the Hungarian National Research Fund, grant no. K60708. WW's research was supported in part by ANR-06-BLAN-00058.
 The cooperation of the authors is facilitated by the French-Hungarian bilateral mobility grant Balaton/02/2008.

\end{document}